\documentclass[12pt, twoside, eqno]{article}
\usepackage{latexsym}
\usepackage{amssymb}
\usepackage{amsfonts}
\begin{document}

\noindent {\bf \large An application of $\Gamma$-semigroups 
techniques to the Green's Theorem }\bigskip

\medskip

\noindent{\bf Niovi Kehayopulu}\bigskip

\bigskip

\bigskip{\small

\noindent {\bf Abstract.} The concept of a $\Gamma$-semigroup has 
been introduced by Mridul Kanti Sen in the Int. Symp., New Delhi, 
1981. It is well known that the Green's relations play an essential 
role in studying the structure of semigroups. In the present paper we 
deal with an application of $\Gamma$-semigroups techniques to the 
Green's Theorem in an attempt to show the way we pass from semigroups 
to $\Gamma$-semigroups. \medskip

\smallskip

}\section{Introduction}  The notion of a $\Gamma$-ring, a 
generalization of the concept of associative rings, has been 
introduced and studied by N. Nobusawa in [7]. $\Gamma$-rings have 
been also studied by W.E. Barnes in [1]. J. Luh studied many 
properties of simple $\Gamma$-rings and primitive $\Gamma$-rings in 
[6]. The concept of a $\Gamma$-semigroup has been introduced by M.K. 
Sen in 1981 as follows: Given two nonempty sets $S$ and $\Gamma$, $S$ 
is called a $\Gamma$-semigroup if the following assertions are 
satisfied:

(1) $a\alpha b\in S$ and $\alpha a\beta\in\Gamma$ and

(2) $(a\alpha b)\beta c=a(\alpha b\beta)c=a\alpha (b\beta c)$\\for 
all $a,b,c\in S$ and all $\alpha,\beta\in\Gamma$ [9]. In 1986 Sen and 
Saha gave a second definition of $\Gamma$-semigroups as follows:
Let $S=\{a,b,c, ......\}$ and $\Gamma=\{\alpha,\beta,\gamma, ...... 
\}$ be two nonempty sets. Then $S$ is called a $\Gamma$-semigroup if

(1) $a\alpha b\in S$ and

(2) $(a\alpha b)\beta c=a\alpha (b\beta c)$\\for all $a,b,c\in S$ and 
all $\alpha,\beta\in\Gamma$ [10]. One can find this definition of 
$\Gamma$-semigroups in [13] where the notion of a radical in 
$\Gamma$-semigroups and the notion of $\Gamma S$-act over a 
$\Gamma$-semigroup have been introduced and in [11] and [12] where 
the notions of regular and orthodox $\Gamma$-semigroups have been 
introduced and studied. But still, we cannot say that $\Gamma$ is a 
set of binary operations on $M$. In [8] N.K. Saha defines the 
$\Gamma$-semigroup as follows: Given two nonempty sets $S$ and 
$\Gamma$, $S$ is called a $\Gamma$-semigroup if there exists a 
mapping $$S\times \Gamma\times S \rightarrow S \mid (a,\gamma,b) 
\rightarrow a\gamma b$$such that $(a\alpha b)\beta c=a\alpha (b\beta 
c)$ for all $a,b,c\in S$ and all $\alpha,\beta\in\Gamma$, and remarks 
that most usual semigroup concepts, in particular regular and inverse 
semigroups, have their analogous for $\Gamma$-semigroups. Defining 
the $\Gamma$-semigroup via mappings, in an expression of the form 
$a_1\gamma_1 a_2\gamma a_2 \, ... \,a_n\gamma_n$, we can put 
parentheses in any place beginning with some $a_i$ and ending in some 
$a_j$. Our aim is to show the way we pass form semigroups to 
$\Gamma$-semigroups and, as the Green's Theorem is a fundamental 
result on semigroups, we chose Green's Theorem for our investigation. 
Green's equivalences for $\Gamma$-semigroups have been studied by 
T.K. Dutta and T.K. Chatterjee in [4], where Green's Lemma for a 
$\Gamma$-semigroup has been also given. Green's Lemma and Green's 
theorem with a different point of view has been given in [2]. Our 
methodology and the aim completely differs from that one in [2], we 
should only mention that for mappings which are mutually inverses we 
do not need to provide independent proof that they are onto. Mutually 
inverse mappings are both (1--1) and onto mappings. The results of 
the present paper are based on the Lemma 2.2 (Green's lemma) and the 
Theorem 2.3 given by Clifford-Preston in [3], and our aim is to show 
the way we pass from the results on semigroups to the corresponding 
results on $\Gamma$-semigroups in general.\section{Main result} For 
two nonempty sets $M$ and $\Gamma$, define $M\Gamma M$ as the set of 
all elements of the form $a\gamma b$, where $a,b\in M$ and 
$\gamma\in\Gamma$. That is, $M\Gamma M:=\{a\gamma b \mid a,b\in M, 
\,\gamma\in\Gamma\}$. Then $M$ is called a $\Gamma$-semigroup [10, 5] 
if the following assertions are satisfied:

(1) $M\Gamma M\subseteq M$.

(2) If $a,b,c,d\in M$ and $\gamma,\mu\in\Gamma$ such that $a=c$, 
$\gamma=\mu$ and $b=d$, then

\hspace{0,7cm}$a\gamma b=c\mu d$.

(3) $(a\gamma b)\mu c=a\gamma (b\mu c)$ for all $a,b,c\in M$ and all 
$\gamma,\mu\in \Gamma$.\\In other words, a $\Gamma$-semigroup is a 
nonempty set $M$ with a set $\Gamma$ of binary operations on $M$, 
satisfying the ``associativity condition" $(a\gamma b)\mu 
c=a\gamma(b\mu c)$ for all $a,b,c\in M$ and all 
$\gamma,\mu\in\Gamma$. If we have only the conditions (1) and (2), 
then this is the definition of a $\Gamma$-groupoid.

Let us first give an example of a $\Gamma$-semigroup based on the 
definition introduced by M.K. Sen and N.K. Saha in 1986 [10]. We give 
an example of order 3 which is easy to be checked.

Consider the set $M=\{a,b,c\}$, and let $\Gamma=\{\gamma, \mu\}$ be 
the set of two binary operations on $M$ defined by:\bigskip

$\begin{array}{*{20}{l}}
\gamma &\vline& a&\vline& b&\vline& c\\
\hline
a&\vline& a&\vline& b&\vline& c\\
\hline
b&\vline& b&\vline& c&\vline& a\\
\hline
c&\vline& c&\vline& a&\vline& b
\end{array}\,\,\,\,\,\,\,\,\,\,\,\,\,\,\,\,\,\,\,\,\,\,\,\,\,\,\,\,\,\,\,\begin{array}{*{20}{l}}
\mu &\vline& a&\vline& b&\vline& c\\
\hline
a&\vline& b&\vline& c&\vline& a\\
\hline
b&\vline& c&\vline& a&\vline& b\\
\hline
c&\vline& a&\vline& b&\vline& c
\end{array}$\bigskip

\noindent Since $(x\rho y)\omega z=x\rho (y\omega z)$ for every 
$x,y,z\in M$ and every $\rho, \omega\in\Gamma$, $M$ is a 
$\Gamma$-semigroup.

For a $\Gamma$-groupoid $G$ the Green's relations $\cal L$ and $\cal 
R$ are defined  as follows:$${\cal R}:=\{(a,b) \mid a=b \mbox { or }
\exists \;x,y\in G \mbox { and } \gamma,\mu\in\Gamma : a\gamma 
x=b,\;b\mu y=a\},$$$${\cal L}:=\{(a,b) \mid a=b \mbox { or } \exists
\;x,y\in G \mbox { and } \gamma,\mu\in\Gamma : x\gamma a=b, \,y\mu
b=a\}.$$For $a\in G$, we denote by $(a)_{\cal R}$, $(a)_{\cal L}$ the 
${\cal R}$-class of $G$ containing the element $a$ and the ${\cal 
L}$-class of $G$ containing the element $a$, respectively. One can 
easily prove that the relation ${\cal R}$ is a left congruence on 
$G$, that is an equivalence relation on $G$ such that $(a,b)\in\cal 
R$ implies $(c\gamma a,c\gamma b)\in\cal R$ for every $c\in G$ and 
every $\gamma\in\Gamma$, and that the relation $\cal L$ is a right 
congruence on $G$, that is, an equivalence relation on $G$ such that 
$(a,b)\in\cal L$ implies $(a\gamma c,b\gamma c)\in\cal L$ for every 
$c\in G$ and every $\gamma\in\Gamma$. By a right (resp. left) ideal 
of $G$ we mean a nonempty subset $A$ of $G$ such that $A\Gamma 
G\subseteq A$ (resp. $G\Gamma A\subseteq A$). If $R(a)$ (resp. 
$L(a)$) denotes the right (resp. left) ideal of $G$ generated by the 
element $a$ of $G$, then we have $R(a)=\{a\}\cup \{a\}\Gamma G$, 
$L(a)=\{a\}\cup G\Gamma\{a\}$, also $a{\cal R}b$ if and only if 
$R(a)=R(b)$ and $a{\cal L}b$ if and only if $L(a)=L(b)$ [4]. Let us 
begin with the Green's Lemma for $\Gamma$-semigroups.\\We first give 
the Lemma 2.2 (Green) as it is stated in [3]. \medskip

\noindent{\bf Lemma 2.2 (Green)} [3] {\it Let a and b be ${\cal 
R}$-equivalent elements of a semigroup S, and let s and $s'$ be the 
elements of $S^1$ such that $as=b$ and $bs'=a$.} ({\it Such elements 
s and $s'$ must exist}). {\it Then the mappings $x\rightarrow xs$ 
$(x\in (a)_{\cal L})$ and $y\rightarrow ys'$ $(y\in (b)_{\cal L})$ 
are mutually inverse, ${\cal R}$-class preserving, one-to-one 
mappings of $(a)_{\cal L}$ upon $(b)_{\cal L}$, and of $(b)_{\cal L}$ 
upon $(a)_{\cal L}$, respectively.}\\According to Clifford-Preston  
the ``$\sigma$ is ${\cal R}$-class preserving" means that ``if $x\in 
(a)_{\cal L}$, then $xs{\cal R}x$".\medskip

\noindent{\bf Remark 1.} If $G$ is a $\Gamma$-groupoid, $a,b,s\in G$, 
$\gamma\in\Gamma$ and $a\gamma s=b$, then $(x,a)\in{\cal L} \mbox { 
implies } (x\gamma s,b)\in{\cal L}$, and the mapping $\sigma : 
(a)_{\cal L} \rightarrow (b)_{\cal L} \mid x \rightarrow x\gamma s$ 
is well defined. Indeed, since $\cal L$ is a right congruence on $G$, 
$(x,a)\in{\cal L}$ implies $(x\gamma s,a\gamma s)\in{\cal L}$, so 
$(x\gamma s,b)\in{\cal L}$. If $x\in (a)_{\cal L}$, then 
$(x,a)\in{\cal L}$, thus $(x\gamma s,b)\in{\cal L}$, so $x\gamma s\in 
(b)_{\cal L}$. $\hfill\Box$\\By Remark 1, by symmetry, we have the 
following\medskip

\noindent{\bf Remark 2.} If $G$ is a $\Gamma$-groupoid, $b,a,s'\in 
G$, $\mu\in\Gamma$ and $b\mu s'=a$, then $(x,b)\in{\cal L}$ implies  
$(x\mu s',a)\in{\cal L}$, and the mapping $\sigma' : (b)_{\cal L} 
\rightarrow (a)_{\cal L} \mid x \rightarrow x\mu s'$ is well defined. 
$\hfill\Box$

If $f$ is a mapping of a set $A$ into the set $B$ and $g$ a mapping 
of $B$ into $A$, $g\circ f$ denotes the mapping of $A$ into $A$ 
defined by $(g\circ f)(x)=g(f(x))$ for every $x\in A$.\medskip

\noindent{\bf Remark 3.} If $G$ is a $\Gamma$-semigroup, $a,b,s,s'\in 
G$ and $\gamma,\mu\in\Gamma$ such that $a\gamma s=b$ and $b\mu s'=a$, 
then the mappings

$\sigma : (a)_{\cal L} \rightarrow (b)_{\cal L} \mid x \rightarrow 
x\gamma s$ and $\sigma' : (b)_{\cal L} \rightarrow (a)_{\cal L} \mid 
x \rightarrow x\mu s'$ \\are mutually inverse mappings. Indeed:

$\sigma'\circ \sigma : (a)_{\cal L} \rightarrow (a)_{\cal L} \mid x 
\rightarrow (x\gamma s)\mu s'$ and $\sigma\circ \sigma' : (b)_{\cal 
L} \rightarrow (b)_{\cal L} \mid y\rightarrow (y\mu s')\gamma s$.\\
$\sigma'\circ \sigma$ is the identity mapping on $(a)_{\cal L}$. 
Indeed: If $d\in (a)_{\cal L}$, then $d=a$ or there exist $t\in G$, 
$\rho\in\Gamma$ such that $d=t\rho a$. If $d=a$, then
$$(\sigma'\circ \sigma)(d)=(\sigma'\circ \sigma)(a)=(a\gamma s)\mu 
s'=b\mu s'=a=d.$$If $t\in G$ and $\rho\in\Gamma$ such that $d=t\rho 
a$, then\begin{eqnarray*}(\sigma'\circ \sigma)(d):&=&(d\gamma s)\mu 
s' =d\gamma (s\mu s')=(t\rho a)\gamma s\mu s'=t\rho (a\gamma s)\mu 
s'\\&=&t\rho b\mu s'=t\rho (b\mu s')=t\rho a=d.\end{eqnarray*}By 
symmetry, $\sigma\circ \sigma'$ is the identity mapping on $(b)_{\cal 
L}$, thus $\sigma$ and $\sigma'$ are mutually inverse mappings. As a 
result, $\sigma$ and $\sigma'$ are (1--1) mappings of $(a)_{\cal L}$ 
onto $(b)_{\cal L}$ and $(b)_{\cal L}$ onto $(a)_{\cal L}$, 
respectively. $\hfill\Box$\medskip

\noindent{\bf Remark 4.} If $G$ is a $\Gamma$-groupoid, $a,b,s,s'\in 
G$, $\gamma,\mu\in\Gamma$,

$\sigma : (a)_{\cal L} \rightarrow (b)_{\cal L} \mid x \rightarrow 
x\gamma s, \;\; \sigma' : (b)_{\cal L} \rightarrow (a)_{\cal L} \mid 
x \rightarrow x\mu s',$\\and the mapping $\sigma'\circ\sigma$ (resp. 
$\sigma\circ\sigma'$) is the identity mapping, then $x\in (a)_{\cal 
L}$ implies $(x\gamma s,x)\in{\cal R},$ and $x\in (b)_{\cal L}$ 
implies $(x\mu s',x)\in{\cal R}.$ Indeed: If $\sigma'\circ\sigma$ is 
the identity mapping and $x\in (a)_{\cal L}$, then $(\sigma'\circ 
\sigma)(x)=x$, so $(x\gamma s)\mu s'=x$. Since $(x\gamma s)\mu s'=x$ 
and $x\gamma s=x\gamma s$, by the definition of $\cal R$, we have 
$(x\gamma s,x)\in\cal R$. The proof of the rest is similar. 
$\hfill\Box$\medskip

The Remarks mentioned above lead to the following lemma which is the 
Green's Lemma for $\Gamma$-semigroups given by T.K. Dutta and T.K. 
Chatterjee in [4]. \medskip

\noindent{\bf Lemma (Green)} [4] {\it Let a and b be $\cal 
R$-equivalent elements of a $\Gamma$-semigroup G, $a\not=b$, $s,s'$ 
the elements of $G$ and $\gamma,\mu$ the elements of $\Gamma$ for 
which $a\gamma s=b$ and $b\mu s'=a$. Then the mappings$$\sigma : 
(a)_{\cal L} \rightarrow (b)_{\cal L} \mid x \rightarrow x\gamma 
s\mbox { and } \sigma' : (b)_{\cal L} \rightarrow (a)_{\cal L} \mid x 
\rightarrow x\mu s'$$are mutually inverse mappings, so $\sigma$ is a 
one-to-one mapping of $(a)_{\cal L}$ onto $(b)_{\cal L}$, and 
$\sigma'$ a one-to-one mapping $(b)_{\cal L}$ onto $(a)_{\cal L}$. 
Moreover, if $x\in (a)_{\cal L}$, then $(x\gamma s,x)\in\cal R$ and 
if $x\in (b)_{\cal L}$, then $(x\mu s',x)\in\cal R$. } 
$\hfill\Box$\medskip

Now let us give the Theorem 2.3 (Green's Theorem) as it is stated in 
[3] :\medskip

\noindent{\bf Theorem 2.3 (Green)} [3] {\it Let a and c be $\cal 
R\circ \cal L$ equivalent elements of a semigroup S. Then there 
exists b in S such that $a{\cal R}b$ and $b{\cal L}c$, and hence 
as=b, $bs'=a$, tb=c, $t'c=b$ for some $s,s',t,t'\in S^1$. The 
mappings:$$\sigma : (a)_{\cal R}\cap (a)_{\cal L} \rightarrow 
(c)_{\cal R}\cap (c)_{\cal L} \mid x \rightarrow t\rho x\gamma 
s$$and$$\sigma' : (c)_{\cal R}\cap (c)_{\cal L} \rightarrow (a)_{\cal 
R}\cap (a)_{\cal L} \mid x \rightarrow t'\zeta x\mu s'.$$are mutually 
inverse, one-to-one mappings of $(a)_{\cal R}\cap (a)_{\cal L}$ and 
$(c)_{\cal R}\cap (c)_{\cal L}$ upon each other.}\medskip

We denote by ${\cal R}\circ {\cal L}$ the relation on $G\times G$ 
defined by:$${\cal R}\circ {\cal L}=\{(a,c) \mid \exists\; b\in G : 
(a,b)\in{\cal R} \mbox { and } (b,c)\in{\cal L}\}.$$
To shorten the proof of Green's Theorem for $\Gamma$-semigroups, 
exactly as in [3], we use the usual notation $G^1=G\cup \{1\}$, and 
define $a\gamma 1=a$ and $1\gamma a=a$ for every $a\in G$.

Using this notation, the relations $\cal R$ and $\cal L$ defined 
above are clearly written as follows:
$${\cal R}=\{(a,b) \mid
\exists \;x,y\in G^1 \mbox { and } \gamma,\mu\in\Gamma : a\gamma 
x=b,\;b\mu y=a\},$$$${\cal L}=\{(a,b) \mid  \exists
\;x,y\in G^1 \mbox { and } \gamma,\mu\in\Gamma : x\gamma a=b, \,y\mu
b=a\}.$$

In the following, we present the Green's Theorem for 
$\Gamma$-semigroups.\medskip

\noindent{\bf Theorem (Green)} (cf. also [2]) {\it Let a and c be 
$\cal R\circ \cal L$ equivalent elements of a $\Gamma$-semigroup $G$, 
b the element of G for which $(a,b)\in\cal R$ and $(b,c)\in\cal L$, 
$s,s',t,t'$ the elements of $G^1$ and $\gamma,\mu,\rho,\zeta$ the 
elements of $\Gamma$ for which $$a\gamma s=b,\;b\mu s'=a, \;t\rho 
b=c, \;t'\zeta c=b.$$We consider the mappings:$$\sigma : (a)_{\cal 
R}\cap (a)_{\cal L} \rightarrow (c)_{\cal R}\cap (c)_{\cal L} \mid x 
\rightarrow t\rho x\gamma s$$and$$\sigma' : (c)_{\cal R}\cap 
(c)_{\cal L} \rightarrow (a)_{\cal R}\cap (a)_{\cal L} \mid x 
\rightarrow t'\zeta x\mu s'.$$Then the mappings $\sigma$ and 
$\sigma'$ are mutually inverse, $\sigma$ is a one-to-one mapping of 
$(a)_{\cal R}\cap (a)_{\cal L}$ onto $(c)_{\cal R}\cap (c)_{\cal L}$ 
and ${\sigma'}$ is a one-to-one mapping of $(c)_{\cal R}\cap 
(c)_{\cal L}$ onto $(a)_{\cal R}\cap (a)_{\cal L}$.} \medskip

\noindent{\bf Proof.} The mapping $\sigma$ is well defined. Indeed: 
Let $x\in (a)_{\cal R}\cap (a)_{\cal L}$. Since $(x,a)\in {\cal R}$, 
we have $x=a\xi z$ and $a=x\xi'z'$ for some $z,z'\in G^1$ and 
$\xi,\xi'\in\Gamma$. Since $(x,a)\in {\cal L}$, we have $x=w\lambda 
a$ and $a=w'\lambda' x$ for some $w,w'\in G^1$ and 
$\lambda,\lambda'\in\Gamma$. Then we have\begin{eqnarray*}c&=&t\rho 
b=t\rho (a\gamma s)=t\rho (x\xi'z')\gamma s=t\rho(w\lambda 
a)\xi'z'\gamma s=t\rho w\lambda (b\mu s')\xi'z'\gamma s\\&=&t\rho 
w\lambda (a\gamma s)\mu s'\xi'z'\gamma s=t\rho (w\lambda a)\gamma 
s\mu s'\xi'z'\gamma s=t\rho x\gamma s\mu (s'\xi'z'\gamma 
s)\end{eqnarray*}and$$t\rho x\gamma s=t\rho (a\xi z)\gamma s=t\rho 
(b\mu s')\xi z\gamma s=(t\rho b)\mu s'\xi z\gamma s=c\mu (s' \xi 
z\gamma s).$$Since $s'\xi'z'\gamma s, \;s' \xi z\gamma s\in G^1$, 
$\;\mu
\in \Gamma$, $\;t\rho x\gamma s\mu (s'\xi'z'\gamma s)=c$ and $\;c\mu 
(s'\xi'z'\gamma s)=t\rho x\gamma s$, we have $(t\rho x\gamma s, 
c)\in\cal R$, so $t\rho x\gamma s\in (c)_{\cal R}$.\\
We also have\begin{eqnarray*}c&=&t\rho b=t\rho (a\gamma s)=t\rho 
(w'\lambda' x)\gamma s=t\rho w'\lambda' (a\xi z)\gamma s=t\rho 
w'\lambda'(b\mu s')\xi z\gamma s\\&=&t\rho w'\lambda'(t'\zeta c)\mu 
s'\xi z\gamma s=t\rho w'\lambda't'\zeta (t\rho b)\mu s'\xi z\gamma 
s=t\rho w'\lambda't'\zeta t\rho (b\mu s')\xi z\gamma s\\&=&t\rho 
w'\lambda't'\zeta t\rho a\xi z\gamma s=t\rho w'\lambda't'\zeta t\rho 
(a\xi z)\gamma s=(t\rho w'\lambda't')\zeta (t\rho x\gamma 
s)\end{eqnarray*}and
$$t\rho x\gamma s=t\rho (w\lambda a)\gamma s=t\rho w\lambda (a\gamma 
s)=t\rho w\lambda b=t\rho w\lambda (t'\zeta c)=(t\rho w\lambda 
t')\zeta c.$$Since $t\rho w'\lambda't'$, $t\rho w\lambda t'\in G^1$, 
$\zeta \in\Gamma$, $(t\rho w'\lambda't')\zeta (t\rho x\gamma s)=c$ 
and $(t\rho w\lambda t')\zeta c=t\rho x\gamma s$, we have $(t\rho 
x\gamma s,c)\in\cal L$, so $t\rho x\gamma s\in (c)_{\cal L}$. Hence 
we obtain $t\rho x\gamma s\in (c)_{\cal R}\cap (c)_{\cal L}$, and the 
mapping $\sigma$ is well defined.\\
The mapping $\sigma'$ is well defined. Indeed: Let $x\in (c)_{\cal 
R}\cap (c)_{\cal L}$. Since $(x,c)\in\cal R$, we have $x=c\xi u$ and 
$c=x\xi'u'$ for some $u,u'\in G^1$ and $\xi,\xi'\in\Gamma$. Since 
$(x,c)\in\cal L$, we have $x=w\lambda c$ and $c=w'\lambda'x$ for some 
$w,w'\in G^1$ and $\lambda,\lambda'\in\Gamma$. We 
have\begin{eqnarray*}a&=&b\mu s'=(t'\zeta c)\mu s'=t'\zeta 
(x\xi'u')\mu s'=t'\zeta (w\lambda c)\xi'u'\mu s'=t'\zeta w\lambda 
(t\rho b)\xi'u'\mu s'\\&=&t'\zeta w\lambda t\rho (a\gamma s)\xi'u'\mu 
s'=t'\zeta w\lambda t\rho (b\mu s')\gamma s\xi'u'\mu s'=t'\zeta 
w\lambda (t\rho b)\mu s'\gamma s\xi'u'\mu s'\\&=&t'\zeta (w\lambda 
c)\mu s'\gamma s\xi'u'\mu s'=(t'\zeta x\mu s')\gamma (s\xi'u'\mu 
s')\end{eqnarray*}and$$t'\zeta x\mu s'=t'\zeta (c\xi u)\mu 
s'=(t'\zeta c)\xi u\mu s'=b\xi u\mu s'=(a\gamma s)\xi u\mu s'=a\gamma 
(s\xi u\mu s').$$Since $s\xi'u'\mu s', \;s\xi u\mu s'\in G^1$, 
$\gamma\in\Gamma$, $(t'\zeta x\mu s')\gamma (s\xi'u'\mu s')=a$ and 
$a\gamma (s\xi u\mu s')=t'\zeta x\mu s'$, we have $(t'\zeta x\mu 
s',a)\in{\cal R}$, so $t'\zeta x\mu s'\in (a)_{\cal R}$.\\
We also have\begin{eqnarray*}a&=&b\mu s'=(t'\zeta c)\mu s'=t'\zeta 
(w'\lambda' x)\mu s'=t'\zeta w'\lambda' (c\xi u)\mu s'=t'\zeta 
w'\lambda' (t\rho b)\xi u\mu s'\\&=&t'\zeta w'\lambda' t\rho (t'\zeta 
c)\xi u\mu s'=t'\zeta w'\lambda' t\rho t'\zeta (c\xi u)\mu 
s'=(t'\zeta w'\lambda' t)\rho (t'\zeta x \mu 
s')\end{eqnarray*}and$$t'\zeta x \mu s'=t'\zeta (w\lambda c)\mu 
s'=t'\zeta w\lambda (t\rho b)\mu s'=t'\zeta w\lambda t\rho (b\mu 
s')=(t'\zeta w\lambda t)\rho a.$$Since $t'\zeta w'\lambda' t, 
\;t'\zeta w\lambda t\in G^1$, $\rho\in\Gamma$, $(t'\zeta w'\lambda' 
t)\rho (t'\zeta x \mu s')=a$ and $(t'\zeta w\lambda t)\rho a=t'\zeta 
x \mu s'$, we have $(t'\zeta x \mu s',a)\in {\cal L}$, so $t'\zeta x 
\mu s'\in (a)_{\cal L}$. Hence we obtain $t'\zeta x\mu s'\in 
(a)_{\cal R}\cap (a)_{\cal L}$, and the mapping $\sigma'$ is well 
defined.\\The mappings $\sigma$ and $\sigma'$ are mutually inverse 
mappings. In fact:
$$\sigma'\circ \sigma: (a)_{\cal R}\cap (a)_{\cal L} \rightarrow 
(a)_{\cal R}\cap (a)_{\cal L} \mid x \rightarrow (\sigma'\circ 
\sigma)(x):=\sigma'(\sigma (x))=t'\zeta (t\rho x\gamma s)\mu s',$$
$$\sigma\circ \sigma': (c)_{\cal R}\cap (c)_{\cal L} \rightarrow 
(c)_{\cal R}\cap (c)_{\cal L} \mid x \rightarrow (\sigma\circ 
\sigma')(x):=\sigma(\sigma' (x))=t\rho (t'\zeta x\mu s')\gamma s.$$
If $d\in (a)_{\cal R}\cap (a)_{\cal L}$, then $(\sigma'\circ 
\sigma)(d)=d$. Indeed: Since $(d,a)\in\cal R$, we have $d=a\beta u$ 
and $a=d\beta'u'$ for some $u,u'\in G^1$ and $\beta,\beta'\in\Gamma$. 
Since $(d,a)\in\cal L$, we have $d=w\delta a$ and $a=w'\delta'd$ for 
some $w,w'\in G^1$ and $\delta,\delta'\in\Gamma$. Then we 
have\begin{eqnarray*}(\sigma'\circ\sigma)(d)&=&t'\zeta (t\rho d\gamma 
s)\mu s'=t'\zeta t\rho (w\delta a)\gamma s\mu s'=t'\zeta t\rho 
w\delta (a\gamma s)\mu s'\\&=&t'\zeta t\rho w\delta (b\mu s')=t'\zeta 
t\rho w\delta a=t'\zeta t\rho (w\delta a)=t'\zeta t\rho d =t'\zeta 
t\rho (a\beta u)\\&=&t'\zeta t\rho (b\mu s')\beta u=t'\zeta (t\rho 
b)\mu s'\beta u=t'\zeta c\mu s'\beta u=(t'\zeta c)\mu s'\beta u\\&=& 
b\mu s'\beta u=(b\mu s')\beta u=a\beta u=d. \end{eqnarray*}If $d\in 
(c)_{\cal R}\cap (c)_{\cal L}$, then $(\sigma\circ \sigma')(d)=d$. 
Indeed: Since $(d,c)\in\cal R$, we have $d=c\beta u$ and 
$c=d\beta'u'$ for some $u,u'\in G^1$ and $\beta,\beta'\in\Gamma$. 
Since $(d,c)\in\cal L$, we have $d=w\delta c$ and $c=w'\delta'd$ for 
some $w,w'\in G^1$ and $\delta,\delta'\in\Gamma$. Then we 
have\begin{eqnarray*}(\sigma\circ \sigma')(d)&=&t\rho (t'\zeta d\mu 
s')\gamma s=t\rho t'\zeta (w\delta c)\mu s'\gamma s=t\rho t'\zeta 
w\delta (t\rho b)\mu s'\gamma s\\&=&t\rho t'\zeta w\delta t\rho (b\mu 
s')\gamma s=t\rho t'\zeta w\delta t\rho a\gamma s=t\rho t'\zeta 
w\delta t\rho (a\gamma s)\\&=&t\rho t'\zeta w\delta t\rho b=t\rho 
t'\zeta w\delta (t\rho b)=t\rho t'\zeta w\delta c=t\rho t'\zeta 
(w\delta c)\\&=&t\rho t'\zeta d=t\rho t'\zeta (c\beta u)=t\rho 
(t'\zeta c)\beta u=t\rho b\beta u\\&=&(t\rho b)\beta u=c\beta 
u=d.\end{eqnarray*}Since the mappings $\sigma$ and $\sigma'$ are 
mutually inverse mappings, the mapping $\sigma$ is an one-to-one 
mapping of $(a)_{\cal R}\cap (a)_{\cal L}$ onto $(c)_{\cal R}\cap 
(c)_{\cal L}$, and the mapping $\sigma'$ is an one-to-one mapping of 
$(c)_{\cal R}\cap (c)_{\cal L}$ onto $(a)_{\cal R}\cap (a)_{\cal L}$. 
$\hfill\Box$ \medskip

With my best thanks to Dr. Michael Tsingelis who helped me in the 
natural and easy proof of the Green's Theorem on semigroups on which 
the present article is based.{\small
\smallskip

\noindent University of Athens, Department of Mathematics, 15784 
Panepistimiopolis, Greece\\email: nkehayop@math.uoa.gr


\begin{thebibliography}{99}
\bibitem{1} W.E. Barnes, On the $\Gamma$-rings of Nobusawa, Pacific 
J. Math. {\bf 18} (1966), 411--422.
\bibitem{2}R. Chinram, P. Siammai, Green's lemma and Green's theorem 
for $\Gamma$-semigroups, Lobachevskii J. Math. 30, no. 3 (2009), 
208--213.
\bibitem{3} A.H. Clifford, G.B. Preston, The Algebraic Theory of 
Semigroups, Vol. I. Mathematical Surveys, No. 7 American 
Mathematical Society, Providence, R.I. 1961.
\bibitem{4} T.K. Dutta, T.K. Chatterjee, Green's equivalences on 
$\Gamma$-semigroup. Bull. Calcutta Math. Soc.  80, no. 1 (1988), 
30--35.
\bibitem{5} N. Kehayopulu, On prime, weakly prime ideals in 
$po$-$\Gamma$-semigroups, Lobachevskii J. Math. {\bf 30}, no. 4 
(2009), 257--262.
\bibitem{6} Jiang Luh, On the theory of simple $\Gamma$-rings, 
Michigan Math. J. {\bf 16} (1969), 65--75.
\bibitem{7} N. Nobusawa, On a generalization of the ring theory, 
Osaka J. Math. {\bf 1} (1964), 81--89.
\bibitem{8} N.K. Saha, The maximum idempotent separating congruence 
on an inverse $\Gamma$-semigroup, Kyungpook Math. J. {\bf 34}, 
no. 1 (1994), 59--66.
\bibitem{9} M.K. Sen, On $\Gamma$-semigroup. In: Algebra and its 
applications, Int. Symp. New Delhi, 1981. Lecture Notes in Pure 
and Appl. Math. {\bf 91}, Dekker, New York, 1984, pp.  301--308.
\bibitem{10} M.K. Sen, N.K. Saha, On $\Gamma$-semigroup. I. Bull. 
Calcutta Math. Soc. {\bf 78}, no. 3 (1986), 180--186.
\bibitem{11} M.K. Sen, N.K. Saha, The maximum idempotent-separating 
congruence on an orthodox $\Gamma$-semigroup, J. Pure Math. {\bf 
7} (1990), 39--47.
\bibitem{12} M.K. Sen, N.K. Saha, Orthodox $\Gamma$-semigroups, 
Internat. J. Math. Math. Sci. {\bf 13}, no. 3 (1990), 527--534.
\bibitem{13} M.K. Sen, A. Seth, Radical of $\Gamma$-semigroup, Bull. 
Calcutta Math. Soc. {\bf 80}, no. 3 (1988), 189--196.
\end{thebibliography}
\end{document}